\documentclass{amsart}
\usepackage{amssymb}
\usepackage{orcidlink}

\newtheorem{theorem}{Theorem}[section]
\newtheorem{lemma}[theorem]{Lemma}
\newtheorem{proposition}[theorem]{Proposition}
\newtheorem{corollary}[theorem]{Corollary}
\theoremstyle{definition}
\newtheorem{definition}[theorem]{Definition}
\newtheorem{example}[theorem]{Example}
\theoremstyle{remark}
\newtheorem{remark}[theorem]{Remark}

\newcommand{\K}{\mathbb{K}}
\newcommand{\R}{\mathbb{R}}
\newcommand{\N}{\mathbb{N}}
\newcommand{\Rm}{\mathbb{R}_{-}}
\newcommand{\sun}{\odot}
\newcommand{\cP}{\mathcal{P}}
\newcommand{\cV}{\mathcal{V}}
\newcommand{\cC}{\mathcal{C}}
\newcommand{\cD}{\mathcal{D}}
\newcommand{\BUC}{\mathrm{BUC}}
\newcommand{\one}{\mathbf{1}}
\newcommand{\dual}[2]{\langle #1,#2\rangle}
\DeclareMathOperator{\co}{co}
\DeclareMathOperator{\spn}{span}
\newcommand{\cco}[1]{\overline{\co}_{#1}\,}

\begin{document}

\title{Operator generated topologies and the Krein property}

\author{Christian Budde \orcidlink{0000-0001-5914-3266}}
\address{University of the Free State, Department of Mathematics and Applied Mathematics, P.O.~Box 339, 9300 Bloemfontein, South Africa}
\address{DSTI-NRF Centre of Excellence in Mathematical and Statistical Sciences (CoE-MaSS), South Africa}
\email{BuddeCJ@ufs.ac.za}

\author{Josef Kreulich \orcidlink{0000-0002-7155-2828}}
\address{Universit\"at Duisburg-Essen}
\email{Josef.Kreulich@gmail.com}

\subjclass[2020]{Primary 46A50; Secondary 46A55, 46B50, 47D06, 34K06}
\keywords{Krein's theorem, convex compactness property, barycentre, norming subspace,
Mackey topology, adjoint semigroup, sun dual, delay equation}

\begin{abstract}
Let $X$ be a Banach space and $Z\subseteq X^*$ a subspace separating the points of $X$.
We say that $Z$ has the Krein property if the $\sigma(X,Z)$-closed convex hull of every
norm-bounded $\sigma(X,Z)$-compact set is $\sigma(X,Z)$-compact. We show that this property
passes from $Z$ to every subspace between $Z$ and its sequential weak$^*$ closure in $X^*$,
and that the closed convex hulls do not change. In particular, it passes to the smallest
$V^*$-invariant space $Z[V]\supseteq Z$ whenever the operator $V$ is a pointwise
$\sigma(X,Z)$-limit of a sequence of $\sigma(X,Z)$-continuous operators; for nets this fails in
general and holds under a double limit condition. As a source of spaces with the Krein property we show
that the Mackey topologies $\mu(X,Z)$ and $\mu(Z,X)$ are complete if the identity of $X$ is
suitably approximated by a sequence of operators whose adjoints map $X^*$ into $Z$; in
particular, $\mu(E^\sun,E)$ and $\mu(E,E^\sun)$ are complete for every $C_0$-semigroup on a
Banach space $E$ with sun dual $E^\sun$. The results are applied to the
history space $\BUC(\Rm)$ of delay equations with infinite delay, where every space of
bounded measures containing $L^1(\Rm)$ has the Krein property.
\end{abstract}

\maketitle

\section{Introduction}

By the theorem of Krein and \v{S}mulian, the closed convex hull of a weakly compact subset of
a Banach space is weakly compact. For a dual pair $\dual{X}{Z}$ the corresponding statement
for $\sigma(X,Z)$ holds if the Mackey topology $\mu(X,Z)$ is quasi-complete
\cite[$\S$24.5(4)]{Koethe}, and it may fail otherwise. We consider a Banach space $X$ and a
subspace $Z\subseteq X^*$ which separates the points of $X$, but need be neither $X^*$ nor a
predual of $X$. Examples are the norming dual pairs of transition semigroups
\cite{Kunze,GerlachKunze}, the pairings induced by adjoint semigroups \cite{vanNeerven}, and
those behind compactifications of semigroup representations \cite{BuddeKreulich}; in each
case compact convex hulls are needed to integrate \cite{Kunze,AndroulakisZiemke} and to
form means.

\begin{definition}\label{def:krein}
A subspace $Z\subseteq X^*$ separating the points of $X$ has the \emph{Krein property} if for
every norm-bounded, $\sigma(X,Z)$-compact set $K\subseteq X$ the $\sigma(X,Z)$-closed convex
hull of $K$ is $\sigma(X,Z)$-compact.
\end{definition}

This is the convex compactness property of \cite{OstlingWilansky} for the norm-bounded
subsets of $(X,\sigma(X,Z))$, and the Krein--Smulian property of the pair $(X,Z)$ in
\cite{AndroulakisZiemke}. Recall that $Z$ is \emph{norming} if there is $c\ge1$ with
$\|x\|\le c\,\sup\{|f(x)|:f\in Z,\ \|f\|\le1\}$ for all $x\in X$, and \emph{$1$-norming} if
one can take $c=1$. The Krein property is a restriction: it fails for $X=\ell^1[0,1]$ and
the norm-closed $1$-norming subspace $Z=C[0,1]$ of $X^*=\ell^\infty[0,1]$
\cite{BonetCascales}. It holds for every $1$-norming $Z$ if $X$ contains no isomorphic copy
of $\ell^1[0,1]$ \cite{CMV,CascalesShvydkoy}, and for every $Z$ if $X$ is weakly compactly
generated \cite{CascalesShvydkoy}; the Banach spaces with $\mu(X,Z)$ complete for all
norm-closed norming $Z$ are studied in \cite{GMZ,GMR}.

In this note we study the stability of the Krein property under refinements of
$\sigma(X,Z)$ that are generated by operators. In the situations above one meets operators
$V$ on $X$ which are pointwise limits of $\sigma(X,Z)$-continuous operators without being
$\sigma(X,Z)$-continuous themselves: ergodic projections, idempotents of a compactification,
solution operators. The coarsest topology of the form $\sigma(X,N)$,
$Z\subseteq N\subseteq X^*$, for which such a $V$ is continuous is given by
$N=Z[V]:=\spn\bigcup_{n\ge 0}V^{*n}Z$. Our main result, Theorems~\ref{thm:seq}
and~\ref{thm:op}, states that the Krein property passes from $Z$ to $Z[V]$ if $V$ is the
pointwise $\sigma(X,Z)$-limit of a \emph{sequence} of $\sigma(X,Z)$-continuous operators,
and that the closed convex hulls of norm-bounded $\sigma(X,Z[V])$-compact sets are the same
for both topologies, which coincide on them. The same holds for sets of such operators and
for every subspace between $Z$ and its sequential weak$^*$ closure $Z_\sigma$. Sequences
cannot be replaced by bounded nets
(Example~\ref{ex:nets}); for nets the result holds under a double limit condition
(Theorem~\ref{thm:nets}, Remark~\ref{rem:nets}).

Section~\ref{sec:sun} provides spaces with the Krein property to start from. Let $Z$ be
norm-closed, and let $(T_n)$ be a bounded sequence in $L(X)$ with $T_n^*X^*\subseteq Z$.
We show that $(X,\mu(X,Z))$ is complete if $T_nx\to x$ for all $x\in X$, and that
$(Z,\mu(Z,X))$ is complete if $T_n^*z\to z$ for all $z\in Z$
(Theorem~\ref{thm:approx}); by Krein's theorem, this yields the Krein property. For a
$C_0$-semigroup on a Banach space $E$ with sun dual $E^\sun\subseteq E^*$ both conditions
hold for means of the semigroup, so that the Mackey topologies $\mu(E^\sun,E)$ on $E^\sun$
and $\mu(E,E^\sun)$ on $E$ are complete (Corollary~\ref{cor:sun},
Remark~\ref{rem:gmr}).

In Section~\ref{sec:delay} we apply the results to the history space $X=\BUC(\Rm)$ of delay
equations with infinite delay. It is the sun dual of $L^1(\Rm)$ for the translation
semigroup, but point delays, the evaluation at the present time and solution operators are
not $\sigma(X,L^1(\Rm))$-continuous. We show that every space $N$ of bounded measures with
$L^1(\Rm)\subseteq N\subseteq M(\Rm)$ has the Krein property, and that every bounded operator
on $X$ whose adjoint maps $L^1(\Rm)$ into $M(\Rm)$ is a limit of the kind considered above
(Theorem~\ref{thm:delay}).

\section{The barycentric criterion}\label{sec:bary}

Throughout, $X$ is a Banach space over $\K\in\{\R,\mathbb{C}\}$ with closed unit ball $B_X$,
and $L(X)$ is the space of bounded linear operators on $X$. All subspaces $G\subseteq X^*$
are assumed to separate the points of $X$, so that $\sigma(X,G)$ is a Hausdorff locally
convex topology with dual $G$; we write $\cco{G}K$ for the $\sigma(X,G)$-closed convex hull
of a set $K\subseteq X$. A linear map $V\colon X\to X$ is $\sigma(X,G)$-continuous if and
only if $g\circ V\in G$ for all $g\in G$. For a compact Hausdorff space $K$ let $\cP(K)$ be
the set of Radon probability measures on $K$ with the weak$^*$ topology of $C(K)^*$; it is
compact and convex, and the finitely supported measures are dense in it.

\begin{lemma}[{\cite[Prop.~1.1]{Phelps}}]\label{lem:bary}
Let $E$ be a Hausdorff locally convex space, $H\subseteq E$ compact and convex,
$K\subseteq H$ compact, and $\nu\in\cP(K)$. Then there is exactly one point $r(\nu)\in H$,
the \emph{barycentre} of $\nu$, with
\[
  \ell(r(\nu))=\int_K \ell\,d\nu\qquad\text{for all }\ell\in E'.
\]
\end{lemma}

For $\K=\mathbb{C}$ one applies \cite{Phelps} to the underlying real space, whose dual
consists of the real parts of the elements of $E'$.

\begin{proposition}[Barycentric criterion; {\cite[Ch.~14, Prop.~14.5, pp.~113--115]{Voigt}}]\label{prop:bary}
Let $Z\subseteq N$ be subspaces of $X^*$, let $K\subseteq X$ be $\sigma(X,Z)$-compact such
that $H:=\cco{Z}K$ is $\sigma(X,Z)$-compact, and let $r(\nu)\in H$, $\nu\in\cP(K)$, be the
barycentres of Lemma~\ref{lem:bary} for $E=(X,\sigma(X,Z))$. Then
$H=\{r(\nu):\nu\in\cP(K)\}$ \cite[Section~1]{Phelps}, and the following assertions are
equivalent:
\begin{itemize}
\item[(a)] $\cco{N}K$ is $\sigma(X,N)$-compact;
\item[(b)] $H$ is $\sigma(X,N)$-compact;
\item[(c)] every $n\in N$ is $\sigma(X,Z)$-continuous on $K$, and
$n(r(\nu))=\int_K n\,d\nu$ for all $n\in N$ and $\nu\in\cP(K)$.
\end{itemize}
If they hold, then $\cco{N}K=H$, and $\sigma(X,N)$ and $\sigma(X,Z)$ coincide on $H$.
\end{proposition}

Indeed, (a)$\Leftrightarrow$(c) is the cited result for $(X,\sigma(X,N))$ and the identity
of $K$, whose Pettis integral with respect to $\nu$ is $r(\nu)$; the rest holds because a
compact topology coincides with every coarser Hausdorff topology.

\section{Sequential closures and operator generated topologies}\label{sec:main}

A subspace $Y\subseteq X^*$ is called \emph{weak$^*$ sequentially closed} if $g\in Y$
whenever $g_k\in Y$ and $g_k(x)\to g(x)$ for all $x\in X$. By the uniform boundedness
principle, such a sequence is norm-bounded and its pointwise limit belongs to $X^*$. For a
subspace $Z\subseteq X^*$ we denote by $Z_\sigma$ the smallest weak$^*$ sequentially closed
subspace of $X^*$ containing $Z$, the \emph{sequential weak$^*$ closure} of $Z$. It is
norm-closed.

\begin{lemma}\label{lem:seq}
Let $Z\subseteq X^*$ be a subspace, and let $K\subseteq X$ be norm-bounded and
$\sigma(X,Z)$-compact such that $H:=\cco{Z}K$ is $\sigma(X,Z)$-compact, with barycentres
$r(\nu)\in H$, $\nu\in\cP(K)$. Then for every $g\in Z_\sigma$ the function $g|_K$ is bounded
and Borel measurable on $(K,\sigma(X,Z))$, and
\begin{equation}\label{eq:bary}
  g(r(\nu))=\int_K g\,d\nu\qquad\text{for all }\nu\in\cP(K).
\end{equation}
\end{lemma}

\begin{proof}
Let $\cC$ be the set of all $g\in X^*$ with these properties. It is a subspace of $X^*$
containing $Z$. Let $g_k\in\cC$ converge pointwise on $X$ to $g$, and put
$C:=\sup_k\|g_k\|$ and $R:=\sup_{x\in K}\|x\|$. Then $g|_K$ is Borel measurable as a
pointwise limit of Borel functions, $|g_k|\le CR$ on $K$, and by dominated convergence
\[
  g(r(\nu))=\lim_{k\to\infty}g_k(r(\nu))=\lim_{k\to\infty}\int_K g_k\,d\nu=\int_K g\,d\nu .
\]
So $\cC$ is weak$^*$ sequentially closed, and therefore $Z_\sigma\subseteq\cC$.
\end{proof}

\begin{theorem}\label{thm:seq}
Let $Z\subseteq X^*$ have the Krein property, and let $N$ be a subspace with
$Z\subseteq N\subseteq Z_\sigma$. Then $N$ has the Krein property. More precisely, if
$K\subseteq X$ is norm-bounded and $\sigma(X,N)$-compact and $H:=\cco{Z}K$, then:
\begin{itemize}
\item[(i)] $\cco{N}K=H$;
\item[(ii)] $H$ is $\sigma(X,N)$-compact, and $\sigma(X,N)$ and $\sigma(X,Z)$ coincide on
$H$;
\item[(iii)] $H=\{r(\nu):\nu\in\cP(K)\}$, where $r(\nu)$ is the unique point of $X$ with
$n(r(\nu))=\int_K n\,d\nu$ for all $n\in N$.
\end{itemize}
\end{theorem}

\begin{proof}
The identity from $(K,\sigma(X,N))$ onto $(K,\sigma(X,Z))$ is a continuous bijection of a
compact space onto a Hausdorff space, hence a homeomorphism. So $K$ is
$\sigma(X,Z)$-compact, every $n\in N$ is $\sigma(X,Z)$-continuous on $K$, and $H$ is
$\sigma(X,Z)$-compact by the Krein property of $Z$. By Lemma~\ref{lem:seq},
$n(r(\nu))=\int_K n\,d\nu$ for all $n\in N$ and $\nu\in\cP(K)$. Thus condition (c) of
Proposition~\ref{prop:bary} holds, and the proposition yields (i)--(iii); the point
$r(\nu)$ is unique because $N$ separates the points of $X$.
\end{proof}

\begin{remark}\label{rem:void}
Theorem~\ref{thm:seq} has content only if $Z\ne Z_\sigma$. If $X=Y^*$ is a dual space and
$Z=Y$ is weakly sequentially complete, for instance $Y=L^1(\mu)$ or the predual of a von
Neumann algebra, then $Z_\sigma=Z$.
\end{remark}

We turn to operators. Let $Z\subseteq X^*$ be a subspace. We say that a linear map
$V\colon X\to X$ satisfies condition (O) with respect to $Z$ if
\begin{itemize}
\item[(O)] there are $\sigma(X,Z)$-continuous linear maps $V_k\colon X\to X$, $k\in\N$,
such that $f(V_kx)\to f(Vx)$ for all $x\in X$ and $f\in Z$.
\end{itemize}
If $\tau$ is any Hausdorff locally convex topology on $X$ with $(X,\tau)'=Z$, then every
pointwise $\tau$-limit of a sequence of $\tau$-continuous linear operators satisfies (O),
since $\tau$-continuous linear operators are $\sigma(X,Z)$-continuous. No boundedness is
assumed in (O):

\begin{lemma}\label{lem:auto}
Let $V$ satisfy \textup{(O)} with respect to $Z$. Then $V\in L(X)$ and
$V^*Z_\sigma\subseteq Z_\sigma$.
\end{lemma}

\begin{proof}
Let $(V_k)$ be as in (O), and let $f\in Z$. As $f\circ V_k\in Z\subseteq X^*$ for all
$f\in Z$ and $Z$ separates the points of $X$, the graph of $V_k$ is closed; so
$V_k\in L(X)$ and $V_k^*f\in Z$. The sequence $(V_k^*f)_k$ converges pointwise on $X$ to
$f\circ V$. Hence $f\circ V\in X^*$, and $f\circ V\in Z_\sigma$. As before, the closed
graph theorem gives $V\in L(X)$. Finally, $\cD:=\{g\in X^*:V^*g\in Z_\sigma\}$ is a
subspace of $X^*$ containing $Z$, and it is weak$^*$ sequentially closed: if $g_k\in\cD$
converge pointwise to $g$, then $V^*g_k\in Z_\sigma$ converge pointwise to $V^*g$. So
$Z_\sigma\subseteq\cD$.
\end{proof}

For a set $\cV$ of linear maps on $X$ satisfying (O) with respect to $Z$ let
\[
  Z[\cV]:=\spn\{f\circ V_1\circ\dots\circ V_n: n\ge 0,\ f\in Z,\ V_1,\dots,V_n\in\cV\},
\]
and $Z[V]:=Z[\{V\}]=\spn\bigcup_{n\ge0}V^{*n}Z$ for a single operator.

\begin{theorem}\label{thm:op}
Let $Z\subseteq X^*$ be a subspace, and let $\cV$ be a set of linear maps on $X$ satisfying
\textup{(O)} with respect to $Z$.
\begin{itemize}
\item[(a)] $Z[\cV]$ is the smallest subspace $N$ of $X^*$ with $Z\subseteq N$ and
$V^*N\subseteq N$ for all $V\in\cV$. Thus $\sigma(X,Z[\cV])$ is the coarsest topology of
the form $\sigma(X,N)$, $Z\subseteq N\subseteq X^*$, for which all $V\in\cV$ are
continuous.
\item[(b)] $Z[\cV]\subseteq Z_\sigma$.
\item[(c)] If $Z$ has the Krein property, then so has every subspace $N$ with
$Z\subseteq N\subseteq\overline{Z[\cV]}$, in particular $Z[\cV]$ and its norm closure
$\overline{Z[\cV]}$, and for norm-bounded $\sigma(X,N)$-compact sets $K$ the assertions
\textup{(i)--(iii)} of Theorem~\ref{thm:seq} hold.
\end{itemize}
\end{theorem}

\begin{proof}
(a) The span of the functionals $f\circ V_1\circ\dots\circ V_n$ contains $Z$ and is
invariant under $g\mapsto g\circ V$ for $V\in\cV$, and every subspace with these two
properties contains all these functionals. The second statement follows, since $V$ is
$\sigma(X,N)$-continuous if and only if $V^*N\subseteq N$.

(b) By Lemma~\ref{lem:auto}, $Z_\sigma$ is invariant under $V^*$ for all $V\in\cV$; now
use (a).

(c) Since $Z_\sigma$ is norm-closed, $\overline{Z[\cV]}\subseteq Z_\sigma$ by (b), and
Theorem~\ref{thm:seq} applies.
\end{proof}

\begin{remark}\label{rem:proj}
If $V$ satisfies (O) and $V^2=V$, then $Z[V]=Z+V^*Z$. By Theorem~\ref{thm:op}, $V$ is
then a continuous projection for $\sigma(X,N)$, where $N$ is $Z+V^*Z$ or its norm closure,
and $N$ has the Krein property if $Z$ has it; cf.\ \cite{BuddeKreulich}.
\end{remark}

The sequence in (O) cannot be replaced by a net, even if the net is bounded and converges
pointwise in norm.

\begin{example}\label{ex:nets}
Let $X=M[0,1]$ be the space of bounded Borel measures on $[0,1]$ and $Z=C[0,1]\subseteq
X^*$, so that $\sigma(X,Z)$ is the weak$^*$ topology of $X=C[0,1]^*$. By Alaoglu's theorem,
$Z$ has the Krein property. Let $V\mu:=\sum_{t}\mu(\{t\})\delta_t$ be the projection onto
the atomic part.

For finite sets $\emptyset\ne\Phi\subseteq[0,1]$ and $\Gamma\subseteq X$ and for $j\in\N$
choose continuous functions $b_t\colon[0,1]\to[0,1]$, $t\in\Phi$, with $b_t(t)=1$ and
pairwise disjoint supports so small that
$|\mu|(\operatorname{supp}b_t\setminus\{t\})\le 1/(j\,|\Phi|)$ for all $\mu\in\Gamma$, and
put $V_\alpha\mu:=\sum_{t\in\Phi}\bigl(\int b_t\,d\mu\bigr)\delta_t$ for
$\alpha=(\Phi,\Gamma,j)$. Then $V_\alpha^*f=\sum_{t\in\Phi}f(t)b_t\in Z$ for $f\in Z$, so
$V_\alpha$ is weak$^*$ continuous, and $\|V_\alpha\|\le1$. Directing the indices by
inclusion in the first two components and by the order of $\N$ in the third, we obtain for
$\mu\in\Gamma$
\[
  \|V_\alpha\mu-V\mu\|\le\frac1j+\sum_{t\notin\Phi}|\mu(\{t\})| ,
\]
so $V_\alpha\mu\to V\mu$ in norm for every $\mu\in X$.

Nevertheless, $N:=Z+V^*Z$, the smallest $V^*$-invariant subspace of $X^*$ containing
$Z$, does not have the Krein property. Let $K:=\{\delta_t:t\in[0,1]\}$. For $f,g\in Z$ we
have $(f+g\circ V)(\delta_t)=f(t)+g(t)$, so $t\mapsto\delta_t$ is continuous into
$(X,\sigma(X,N))$, and $K$ is $\sigma(X,N)$-compact. The Riemann sums
$\mu_k:=\frac1k\sum_{i=1}^k\delta_{i/k}\in\co K$ converge weak$^*$ to the Lebesgue measure
$\lambda$, but the total atomic mass $m:=\one\circ V\in N$ satisfies $m(\mu_k)=1$ and
$m(\lambda)=0$; cf.\ Choquet's example \cite[Section~14]{Phelps}. So $(\mu_k)$ has no
$\sigma(X,N)$-cluster point, and $\cco{N}K$ is not $\sigma(X,N)$-compact. In particular,
$V^*Z\not\subseteq Z_\sigma$, and $V$ does not satisfy (O). One may verify that for
$f=\one$ and $n=0$ the double limit condition \eqref{eq:dlc} of Remark~\ref{rem:nets} does
not hold.
\end{example}

For nets the conclusions of Theorem~\ref{thm:op} are equivalent to a weak convergence
condition, for which a double limit criterion is available.

\begin{theorem}[Nets]\label{thm:nets}
Let $Z\subseteq X^*$ have the Krein property, let $(V_\lambda)_{\lambda\in\Lambda}$ be a
net in $L(X)$ with $V_\lambda^*Z\subseteq Z$ and $C:=\sup_\lambda\|V_\lambda\|<\infty$, and
let $V\colon X\to X$ be linear with $f(V_\lambda x)\to f(Vx)$ for all $x\in X$ and
$f\in Z$. Then $V\in L(X)$. Put $N:=\spn\bigcup_{n\ge0}V^{*n}Z$, and let $K\subseteq X$ be
norm-bounded and $\sigma(X,N)$-compact. Then the assertions \textup{(i)--(iii)} of
Theorem~\ref{thm:seq} hold for $K$ and $N$ if and only if
\begin{equation}\label{eq:weak}
  V^{*n}V_\lambda^*f|_K\longrightarrow V^{*(n+1)}f|_K\qquad\text{weakly in }C(K)
\end{equation}
for all $f\in Z$ and $n\ge0$. If $V^2=V$, then \eqref{eq:weak} is only needed for $n=0$.
\end{theorem}

\begin{proof}
For $f\in Z$ we have $|f(Vx)|\le C\|f\|\,\|x\|$, so $f\circ V\in X^*$, and $V\in L(X)$ by
the closed graph theorem, as in the proof of Lemma~\ref{lem:auto}. As in the proof of
Theorem~\ref{thm:seq}, $K$ is $\sigma(X,Z)$-compact, every element of $N$ is
$\sigma(X,Z)$-continuous on $K$, and $H:=\cco{Z}K$ is $\sigma(X,Z)$-compact, with
barycentres $r(\nu)$. By Proposition~\ref{prop:bary}, \textup{(i)--(iii)} hold if and only
if \eqref{eq:bary} holds for all $g\in V^{*n}Z$, $n\ge0$; for $n=0$ this is the definition
of $r(\nu)$. Let $f\in Z$, $g:=V^{*(n+1)}f$ and
$g_\lambda:=V^{*n}V_\lambda^*f\in V^{*n}Z$, so that $g_\lambda\to g$ pointwise on $X$. If
\eqref{eq:bary} holds on $V^{*n}Z$, then it holds for $g$ if and only if
\[
  \int_Kg\,d\nu=g(r(\nu))=\lim_\lambda g_\lambda(r(\nu))=\lim_\lambda\int_Kg_\lambda\,d\nu
  \qquad\text{for all }\nu\in\cP(K),
\]
that is, if and only if $g_\lambda|_K\to g|_K$ weakly in $C(K)$, because $\cP(K)$ spans
$C(K)^*$. Induction on $n$ gives the equivalence. If $V^2=V$, then $V^{*n}Z=V^*Z$ for
$n\ge1$, and only the step from $n=0$ to $n=1$ is needed.
\end{proof}

\begin{remark}\label{rem:nets}
(a) If the bounded set $A:=\{V^{*n}V_\lambda^*f|_K:\lambda\in\Lambda\}$ is relatively
weakly compact in $C(K)$, then \eqref{eq:weak} holds for $f$ and $n$: the weak closure $L$
of $A$ is compact, hence closed, for the coarser Hausdorff topology of pointwise
convergence on $K$, and the two topologies coincide on $L$. By Grothendieck's double limit
theorem \cite[Th\'eor\`eme~6]{Grothendieck}, $A$ is relatively weakly compact if and only
if, for a dense set $D\subseteq K$,
\begin{equation}\label{eq:dlc}
  \lim_i\lim_j\,\dual{V^{*n}V_{\lambda_i}^*f}{x_j}
  =\lim_j\lim_i\,\dual{V^{*n}V_{\lambda_i}^*f}{x_j}
\end{equation}
for all sequences $(\lambda_i)$ in $\Lambda$ and $(x_j)$ in $D$ for which both iterated
limits exist.

(b) If $\Lambda$ has a cofinal sequence, then $V$ satisfies \textup{(O)} and
Theorem~\ref{thm:op} applies; so Theorem~\ref{thm:nets} is of interest only for nets
without cofinal sequences, as in Example~\ref{ex:nets}.

(c) If $K$ is separable, e.g.\ if $Z$ is norm-separable, then \eqref{eq:weak} already
follows if \eqref{eq:dlc} holds for increasing sequences
$\lambda_1\le\lambda_2\le\cdots$. Indeed, with $g_\lambda$ and $g$ as in the proof of
Theorem~\ref{thm:nets}, suppose that \eqref{eq:weak} fails. Then there are
$\nu\in C(K)^*$, $\varepsilon>0$ and a cofinal set $B\subseteq\Lambda$ with
$|\int_K(g_\lambda-g)\,d\nu|\ge\varepsilon$ for $\lambda\in B$. For a dense sequence
$(x_k)$ in $K$ choose $\lambda_1\le\lambda_2\le\cdots$ in $B$ with
$|g_{\lambda_i}(x_k)-g(x_k)|\le1/i$ for $k\le i$. Every sequence in
$A_0:=\{g_{\lambda_i}|_K:i\in\N\}$ has a constant subsequence or one with increasing
indices, so $A_0$ is relatively weakly compact by \eqref{eq:dlc} and \cite{Grothendieck}.
As in (a), now for pointwise convergence on $\{x_k:k\in\N\}$, we get
$g_{\lambda_i}|_K\to g|_K$ weakly, a contradiction.
\end{remark}

\section{Mackey completeness and the sun dual pairing}\label{sec:sun}

The following form of Krein's theorem reduces the Krein property to completeness.

\begin{theorem}[{\cite[$\S$24.5(4)]{Koethe}}]\label{thm:krein}
Let $\dual{F}{G}$ be a dual pair such that $(F,\mu(F,G))$ is quasi-complete, and let
$K\subseteq F$ be $\sigma(F,G)$-compact. Then the $\sigma(F,G)$-closed absolutely convex
hull of $K$ is $\sigma(F,G)$-compact.
\end{theorem}

Both Mackey topologies of a pair $\dual{X}{Z}$ are complete if the identity of $X$ is
suitably approximated by operators $T_n$ with $T_n^*X^*\subseteq Z$. This is
Theorem~\ref{thm:approx} below; its proof uses the following lemma.

\begin{lemma}\label{lem:pairs}
Let $G$ be a Banach space, $F\subseteq G^*$ a subspace separating the points of $G$, and
$v\colon G\to\K$ a linear functional which is $\sigma(G,F)$-continuous on every absolutely
convex $\sigma(G,F)$-compact subset of $G$. Then $v(x_n)\to0$ for every norm-bounded
sequence $(x_n)$ in $G$ with $x_n\to0$ in $\sigma(G,F)$. In particular, $v\in G^*$.
\end{lemma}

\begin{proof}
The operator $J\colon\ell^1\to G$, $Jc:=\sum_nc_nx_n$, is bounded. For $\varphi\in F$ we
have $\dual{\varphi}{Jc}=\sum_nc_n\dual{\varphi}{x_n}$ with
$(\dual{\varphi}{x_n})_n\in c_0$, so $J$ is continuous from $(\ell^1,\sigma(\ell^1,c_0))$
to $(G,\sigma(G,F))$. By Alaoglu's theorem, $K:=J(B_{\ell^1})$ is $\sigma(G,F)$-compact,
and it is absolutely convex. The unit vectors $u_n$ of $\ell^1$ tend to $0$ in
$\sigma(\ell^1,c_0)$; hence $x_n=Ju_n\to0$ in $(K,\sigma(G,F))$, and $v(x_n)\to v(0)=0$.
If $v$ were unbounded on $B_G$, there would be $y_n\in B_G$ with $|v(y_n)|\ge4^n$; the
sequence $x_n:=2^{-n}y_n$ is norm-null, hence $\sigma(G,F)$-null and bounded, but
$|v(x_n)|\ge2^n$.
\end{proof}

\begin{theorem}\label{thm:approx}
Let $Z\subseteq X^*$ be a norm-closed subspace, and let $(T_n)_{n\in\N}$ be a bounded
sequence in $L(X)$ with $T_n^*X^*\subseteq Z$ for all $n\in\N$.
\begin{itemize}
\item[(a)] If $T_nx\to x$ in norm for every $x\in X$, or if $X$ is weakly sequentially
complete and $T_nx\to x$ weakly for every $x\in X$, then $(X,\mu(X,Z))$ is complete. In
particular, $Z\subseteq X^*$ has the Krein property.
\item[(b)] If $\|T_n^*z-z\|\to0$ for every $z\in Z$, then $(Z,\mu(Z,X))$ is complete. In
particular, $X\subseteq Z^*$ has the Krein property.
\end{itemize}
\end{theorem}

\begin{proof}
By Grothendieck's completeness theorem \cite[$\S$21.9]{Koethe}, for a dual pair
$\dual{F}{G}$ the space $(F,\mu(F,G))$ is complete if and only if every linear functional
$v$ on $G$ which is $\sigma(G,F)$-continuous on all absolutely convex
$\sigma(G,F)$-compact subsets of $G$ is given by an element of $F$. By
Theorem~\ref{thm:krein}, completeness implies the Krein properties.

(a) Let $v$ be such a functional for $F=X$, $G=Z$. As $T_n^*$ is weak$^*$ continuous,
$T_n^*(B_{X^*})$ is an absolutely convex $\sigma(Z,X)$-compact subset of $Z$. So
$v\circ T_n^*$ is a linear functional on $X^*$ whose restriction to $B_{X^*}$, and hence
to every bounded subset of $X^*$, is weak$^*$ continuous. As $X$ is complete for its norm
topology $\mu(X,X^*)$, the theorem just quoted, applied to the pair $F=X$, $G=X^*$, yields
$x_n\in X$ with $x^*(x_n)=v(T_n^*x^*)$ for all $x^*\in X^*$. Let $n_k,m_k\to\infty$ and
$x_k^*\in B_{X^*}$, and put $y_k:=(T_{n_k}^*-T_{m_k}^*)x_k^*\in Z$, so that
$x_k^*(x_{n_k}-x_{m_k})=v(y_k)$. The sequence $(y_k)$ is norm-bounded. If $T_nx\to x$ in
norm for all $x\in X$, then $|y_k(x)|\le\|T_{n_k}x-T_{m_k}x\|\to0$, hence $y_k\to0$ in
$\sigma(Z,X)$, and $v(y_k)\to0$ by Lemma~\ref{lem:pairs}; as $n_k$, $m_k$ and $x_k^*$ were
arbitrary, $(x_n)$ is a Cauchy sequence in norm. If $T_nx\to x$ only weakly, we still have
$y_k\to0$ in $\sigma(Z,X)$ for constant sequences $x_k^*=x^*$; so $(x_n)$ is weakly Cauchy,
and it converges weakly if $X$ is weakly sequentially complete. Let $x\in X$ be the limit.
For $z\in Z$ the sequence $(T_n^*z-z)_n$ is norm-bounded and tends to $0$ in
$\sigma(Z,X)$, so that $z(x)=\lim_nz(x_n)=\lim_nv(T_n^*z)=v(z)$ by
Lemma~\ref{lem:pairs}. Thus $v$ is given by $x$.

(b) Let $v$ be such a functional for $F=Z$, $G=X$. By Lemma~\ref{lem:pairs}, $v\in X^*$.
Suppose that $\|T_n^*v-v\|$ does not tend to $0$. Then there are $\varepsilon>0$,
$n_1<n_2<\cdots$ and $e_k\in B_X$ such that $|v(x_k)|>\varepsilon$ for
$x_k:=T_{n_k}e_k-e_k$. The sequence $(x_k)$ is norm-bounded, and
$|z(x_k)|\le\|T_{n_k}^*z-z\|\to0$ for $z\in Z$. So $x_k\to0$ in $\sigma(X,Z)$, and
$v(x_k)\to0$ by Lemma~\ref{lem:pairs}, a contradiction. Hence $v=\lim_nT_n^*v$ in norm,
and $v\in Z$, because $T_n^*v\in Z$ and $Z$ is norm-closed.
\end{proof}

We apply Theorem~\ref{thm:approx} to adjoint semigroups. Let $(S(t))_{t\ge0}$ be a
$C_0$-semigroup on a Banach space $E$, not necessarily bounded. The \emph{sun dual}
\[
  E^\sun:=\bigl\{\varphi\in E^*:\lim_{t\downarrow0}\|S^*(t)\varphi-\varphi\|=0\bigr\}
\]
is a norm-closed, $S^*(t)$-invariant subspace of $E^*$, on which the adjoints $S^*(t)$
form a $C_0$-semigroup \cite[Ch.~1]{vanNeerven}; the topology $\sigma(E,E^\sun)$ is
studied in \cite[Ch.~2]{vanNeerven}. Every $e\in E$ acts on $E^\sun$ by
$\varphi\mapsto\dual{\varphi}{e}$. For $t>0$ let $C_te:=\frac1t\int_0^tS(s)e\,ds$,
$e\in E$. We use the following standard facts, for which we refer again to
\cite[Ch.~1]{vanNeerven}:
\begin{itemize}
\item[(S1)] $C_t^*E^*\subseteq E^\sun$ for all $t>0$.
\item[(S2)] $\|C_t^*\varphi-\varphi\|\le\sup_{0\le s\le t}\|S^*(s)\varphi-\varphi\|\to0$
as $t\downarrow0$ for every $\varphi\in E^\sun$.
\item[(S3)] There is $M\ge1$ with
$\|e\|\le M\sup\{|\dual{\varphi}{e}|:\varphi\in E^\sun,\ \|\varphi\|\le1\}$ for all
$e\in E$. Hence $E^\sun$ is a norming subspace of $E^*$, and $E$ is a norm-closed subspace
of $(E^\sun)^*$.
\end{itemize}

\begin{corollary}\label{cor:sun}
Let $(S(t))_{t\ge0}$ be a $C_0$-semigroup on a Banach space $E$. Then the spaces
$(E^\sun,\mu(E^\sun,E))$ and $(E,\mu(E,E^\sun))$ are complete. Hence the
$\sigma(E^\sun,E)$-closed absolutely convex hull of every $\sigma(E^\sun,E)$-compact subset
of $E^\sun$ is $\sigma(E^\sun,E)$-compact, and the $\sigma(E,E^\sun)$-closed absolutely
convex hull of every $\sigma(E,E^\sun)$-compact subset of $E$ is $\sigma(E,E^\sun)$-compact.
In particular, $E\subseteq(E^\sun)^*$ and $E^\sun\subseteq E^*$ have the Krein property.
\end{corollary}

\begin{proof}
Let $X:=E$, $Z:=E^\sun$ and $T_n:=C_{1/n}$. Then $\|T_n\|\le\sup_{0\le s\le1}\|S(s)\|$,
$T_n^*E^*\subseteq E^\sun$ by (S1), $T_ne\to e$ in norm for every $e\in E$ by strong
continuity, and $\|T_n^*\varphi-\varphi\|\to0$ for every $\varphi\in E^\sun$ by (S2). So
Theorem~\ref{thm:approx} applies, and the statements on the hulls follow from
Theorem~\ref{thm:krein}.
\end{proof}

\begin{remark}\label{rem:gmr}
(a) Under the hypotheses of Theorem~\ref{thm:approx}(a), every $x^*\in X^*$ is the
pointwise limit of the sequence $(T_n^*x^*)$ in $Z$; hence $Z_\sigma=X^*$. So the Krein
property passes there from $Z_\sigma$, where it holds by the Krein--\v{S}mulian theorem,
down to $Z$, whereas Theorem~\ref{thm:seq} passes from $Z$ up to $Z_\sigma$.

(b) By \cite{GMR}, $(X,\mu(X,Y))$ is complete for every norm-closed norming subspace
$Y\subseteq X^*$ if the Banach space $X$ has Efremov's property~(E), in particular if $X$
is separable or weakly compactly generated. By (S3), this covers
$(E,\mu(E,E^\sun))$ if $E$ has property~(E), and $(E^\sun,\mu(E^\sun,E))$ if $E^\sun$ has
it. So Corollary~\ref{cor:sun} goes beyond \cite{GMR} only if the space carrying the
Mackey topology fails property~(E), as $E^\sun=\BUC(\Rm)$ in Section~\ref{sec:delay} does,
which contains $\ell^1(\omega_1)$ \cite{GMR}. For general norm-closed norming subspaces in
place of $E^\sun$ both statements may fail, as the pair $\dual{\ell^1[0,1]}{C[0,1]}$ of
\cite{BonetCascales} shows when read in either direction.
\end{remark}

\section{Delay functionals on the history space}\label{sec:delay}

Linear delay equations with infinite delay,
\[
  u'(t)=\int_{\Rm}u(t+\theta)\,d\eta(\theta)\quad(t\ge0),\qquad u|_{\Rm}=\varphi,
\]
are studied on spaces of histories such as $X:=\BUC(\Rm)$, the bounded uniformly continuous
$\K$-valued functions on $\Rm:=(-\infty,0]$ with the supremum norm; see
\cite{RuessSummers} for the semigroup approach on such spaces and \cite{DiekmannEtAl} for
the sun-star calculus of delay equations. The right-hand side is given
by a bounded Borel measure $\eta$ on $\Rm$, and point delays correspond to atoms of $\eta$.
We identify $\eta\in M(\Rm)$ with the functional $\varphi\mapsto\int\varphi\,d\eta$ on $X$
and $h\in L^1(\Rm)$ with the measure $h(\theta)\,d\theta$, so that
$L^1(\Rm)\subseteq M(\Rm)\subseteq X^*$.

On $E:=L^1(\Rm)$ let $(S(t)h)(\theta):=h(\theta+t)$ for $\theta\le-t$ and $(S(t)h)(\theta):=0$
for $-t<\theta\le0$. This defines a $C_0$-semigroup of isometries with adjoint
$S^*(t)g=g(\cdot-t)$ on $L^\infty(\Rm)$, and
$(C_t^*g)(\theta)=\frac1t\int_0^tg(\theta-s)\,ds$. We put $Q_k:=C_{1/k}^*|_X$ for $k\in\N$.

\begin{lemma}\label{lem:buc}
\begin{itemize}
\item[(a)] $E^\sun=\BUC(\Rm)$, with the supremum norm.
\item[(b)] $Q_k\in L(X)$, and $\|Q_k\varphi-\varphi\|_\infty\to0$ for every $\varphi\in X$.
\item[(c)] $Q_k^*\eta\in L^1(\Rm)$ for every $\eta\in M(\Rm)$. In particular, $Q_k$ is
$\sigma(X,L^1(\Rm))$-continuous.
\end{itemize}
\end{lemma}

\begin{proof}
(a) $\BUC(\Rm)\subseteq E^\sun$ by uniform continuity. Conversely, for $g\in E^\sun$ each
$C_t^*g$ is bounded and Lipschitz continuous, and $\|C_t^*g-g\|_\infty\to0$ by (S2); so
$g\in\BUC(\Rm)$. Now (b) follows from (a) and (S2). For (c), Fubini's theorem gives
\[
  \int_{\Rm}Q_k\varphi\,d\eta=\int_{\Rm}\varphi(\rho)\,h_k(\rho)\,d\rho\qquad\text{with }
  h_k(\rho):=k\,\eta\bigl([\rho,\rho+\tfrac1k]\cap\Rm\bigr),
\]
and $\|h_k\|_1\le\|\eta\|$.
\end{proof}

By Lemma~\ref{lem:buc}(a) and Corollary~\ref{cor:sun}, $L^1(\Rm)$ has the Krein property as
a space of functionals on $X$. Since $X$ contains $\ell^\infty$, and hence
$\ell^1[0,1]$, the criteria of \cite{CMV,CascalesShvydkoy} quoted in the introduction do
not apply.

\begin{theorem}\label{thm:delay}
Let $X=\BUC(\Rm)$.
\begin{itemize}
\item[(a)] Every subspace $N$ with $L^1(\Rm)\subseteq N\subseteq M(\Rm)$ has the Krein
property, and for norm-bounded $\sigma(X,N)$-compact sets $K\subseteq X$ the assertions
\textup{(i)--(iii)} of Theorem~\ref{thm:seq} hold with $Z=L^1(\Rm)$.
\item[(b)] Let $V\in L(X)$ with $V^*L^1(\Rm)\subseteq M(\Rm)$. Then the operators $VQ_k$
are $\sigma(X,L^1(\Rm))$-continuous, and $VQ_k\varphi\to V\varphi$ in norm for every
$\varphi\in X$. Thus $V$ satisfies \textup{(O)} with respect to $Z=L^1(\Rm)$, and $Z[V]$
has the Krein property. If, moreover, $V^*M(\Rm)\subseteq M(\Rm)$, then
$Z[V]\subseteq M(\Rm)$.
\end{itemize}
\end{theorem}

\begin{proof}
(a) Let $\eta\in M(\Rm)$. By Lemma~\ref{lem:buc}, $Q_k^*\eta\in L^1(\Rm)$ and
$(Q_k^*\eta)(\varphi)=\eta(Q_k\varphi)\to\eta(\varphi)$ for all $\varphi\in X$. Hence
$M(\Rm)\subseteq L^1(\Rm)_\sigma$, and Theorem~\ref{thm:seq} applies.

(b) We have $(VQ_k)^*h=Q_k^*V^*h\in L^1(\Rm)$ for $h\in L^1(\Rm)$ by
Lemma~\ref{lem:buc}(c), and $\|VQ_k\varphi-V\varphi\|\le\|V\|\,\|Q_k\varphi-\varphi\|\to0$.
The rest follows from Theorem~\ref{thm:op}.
\end{proof}

\begin{example}[Delay functionals]\label{ex:rankone}
Let $0\ne\eta\in M(\Rm)$ and $0\ne u\in X$, and let $V\varphi:=\eta(\varphi)\,u$. Then
$V^*g=g(u)\,\eta$ for $g\in X^*$, so Theorem~\ref{thm:delay}(b) applies, and
$Z[V]=L^1(\Rm)+\K\eta$. The operator $V$ is $\sigma(X,L^1(\Rm))$-continuous only if $\eta$
is absolutely continuous. As $V^2=\eta(u)V$, it is a projection only if $\eta(u)=1$. For
$\eta=\delta_0$ the space $Z[V]=L^1(\Rm)\oplus\K\delta_0$ adjoins the evaluation at the
present time, for $\eta=\delta_{-r}$ a point delay. By Theorem~\ref{thm:delay}(a), any set
of delay functionals can be adjoined to $L^1(\Rm)$ without losing the Krein property.
\end{example}

\begin{example}[A delay semigroup]\label{ex:semigroup}
Let $a\in\K$ and $r>0$, and let $T(t)\varphi:=u(t+\cdot)|_{\Rm}$, where $u$ is the solution
of $u'(t)=a\,u(t-r)$ for $t\ge0$ with $u|_{\Rm}=\varphi\in X$. For $0\le t\le r$,
\[
  (T(t)\varphi)(\theta)=
  \begin{cases}
  \varphi(t+\theta), & \theta\le-t,\\[2pt]
  \varphi(0)+a\displaystyle\int_{-r}^{t+\theta-r}\varphi(\rho)\,d\rho, & -t<\theta\le0,
  \end{cases}
\]
so that $T(t)\in L(X)$, and for $h\in L^1(\Rm)$
\[
  T(t)^*\delta_0=\delta_0+a\,\one_{[-r,t-r]},\qquad
  T(t)^*h=h(\cdot-t)+a\,\kappa+\Bigl(\int_{-t}^0h(\theta)\,d\theta\Bigr)\delta_0,
\]
where $\kappa(\rho):=\one_{[-r,t-r]}(\rho)\int_{\rho+r-t}^0h(\theta)\,d\theta$. Hence,
for $0<t\le r$, $T(t)$ is not $\sigma(X,L^1(\Rm))$-continuous and not a projection, and
$Z[T(t)]=L^1(\Rm)\oplus\K\delta_0$. By the semigroup law, every $T(t)$ is continuous for
the topology generated by $L^1(\Rm)\oplus\K\delta_0$. For $a\ne0$ the delay functional
$a\,\delta_{-r}$ of the equation is not; it is continuous for $\sigma(X,N)$ whenever $\delta_{-r}\in N$. By
Theorem~\ref{thm:delay}(a), all these spaces $N\subseteq M(\Rm)$ have the Krein property.
\end{example}

\begin{remark}\label{rem:compare}
(a) For $a=0$ the adjoints $T(t)^*$ of Example~\ref{ex:semigroup} restrict to a
$C_0$-semigroup on $L^1(\Rm)\oplus\K\delta_0$ with adjoint
$(g,\gamma)\mapsto(g(\cdot+t)\one_{(-\infty,-t]}+\gamma\one_{(-t,0]},\gamma)$ and sun dual
$X$. So for this $N$ the Krein property also follows from Corollary~\ref{cor:sun}. We do
not know such a representation for spaces $N$ containing other point masses.

(b) For $N=M(\Rm)$ and the larger space $C_b(\Rm)$ in place of $X$, the statement
corresponding to Theorem~\ref{thm:delay}(a) is classical: $M(\Rm)$ is the dual of
$C_b(\Rm)$ with the strict topology, which is complete \cite{Buck} and coincides with the
Mackey topology \cite{Conway}, so that Theorem~\ref{thm:krein} applies. However, $X$ is a
proper dense subspace of $C_b(\Rm)$ for the strict topology, and the point of
Theorem~\ref{thm:delay} is that the closed convex hulls do not leave $X$.
\end{remark}

\section*{Declarations}

\subsection*{Funding}
C.B acknowledges research support by the National Research Foundation (Grant number: SRUG220317136) as well as the DSTI-NRF Centre of Excellence in Mathematical and Statistical Sciences (CoE-MaSS), South Africa (Grant number: 2026-019-OPA-Evolution Equations). It is acknowledged that opinions, findings and conclusions or recommendations expressed in any publication generated by this supported research is that of the author(s). The National Research Foundation accepts no liability whatsoever in this regard. Opinions expressed and conclusions arrived at are those of the author and are not necessarily to be attributed to the CoE-MaSS.

\subsection*{Author Contributions} 
The authors reviewed the manuscript.

\subsection*{Conflict of Interest}
There is no conflict of interest.

\end{document}